\documentclass[psamsfonts,reqno]{amsart}

\usepackage{graphicx}
\usepackage{float}
\usepackage{amssymb,amsfonts}
\usepackage[all,arc]{xy}
\usepackage{enumerate}
\usepackage{mathrsfs}
\usepackage[margin=1.1in]{geometry}
\usepackage{xcolor}
\usepackage{url}
\usepackage{pdfpages}
\usepackage{arydshln}
\usepackage{tikz-cd}
\usepackage{fontawesome5}
\usepackage{hyperref}

\definecolor{darkgreen}{rgb}{0.0,0.5,0.0}
\definecolor{darkblue}{rgb}{0.0,0.0,0.5}
\definecolor{darkred}{rgb}{0.65,0.0,0.0}

\hypersetup{
  colorlinks=true,
  linkcolor=darkred,
  citecolor=darkgreen,
  urlcolor=darkblue
}

\theoremstyle{plain}
\newtheorem{theorem}{Theorem}[section]
\newtheorem{lemma}[theorem]{Lemma}
\newtheorem{proposition}[theorem]{Proposition}
\newtheorem{corollary}[theorem]{Corollary}

\theoremstyle{definition}

\theoremstyle{remark}
\newtheorem{remark}[theorem]{Remark}

\numberwithin{equation}{section}

\makeatletter
\renewcommand{\@biblabel}[1]{%
  \textcolor{black}{[}\textcolor{darkgreen}{#1}\textcolor{black}{]}}

\renewenvironment{thebibliography}[1]
  {\section*{\refname}%
   \begin{list}{\@biblabel{\arabic{enumiv}}}{%
      \usecounter{enumiv}%
      \setlength{\labelwidth}{4em}%
      \setlength{\labelsep}{1em}%
      \setlength{\leftmargin}{\labelwidth}%
      \setlength{\rightmargin}{0pt}%
      \setlength{\itemsep}{1em}%
      \setlength{\parsep}{0pt}}}
  {\end{list}}
\makeatother

\title{On Some Applications to the Newton--Raphson Method}

\author{Treanungkur Mal}

\date{}

\begin{document}

\maketitle

\pagenumbering{gobble}

\begin{abstract}
Newton's method is classically viewed as an iterative procedure for
approximating a zero of a nonlinear function. In this paper, we consider
the Newton iterates from a quadrature perspective: the iterates generate a partition of the interval between the initial point and the
root, which in turn induces a composite trapezoidal rule. We establish an
exact decomposition of the integral associated with this Newton-generated
quadrature and characterize the resulting quadrature error. We also obtain
an explicit bound for this error in terms of the initial Newton decrement
and give a condition under which it is smaller than the error of the trapezoidal rule on the original interval. These results provide
a connection between Newton's method and the quadrature rule
induced by its iterates.
\end{abstract}

\tableofcontents

\pagenumbering{arabic}


\section{Introduction}

Newton's method is one of the classical methods for finding approximate
solutions of nonlinear equations. The basic
idea is to replace a nonlinear equation locally by its linear
approximation and use the resulting linear equation to obtain a new
approximation. Given a differentiable function $f$ and an initial guess
$x_0$, this leads to the iteration
\[
x_{n+1}=x_n-\frac{f(x_n)}{f'(x_n)}.
\]
The convergence properties of this method have been studied extensively,
and under suitable hypotheses Newton's method exhibits quadratic
convergence to a simple root; see, for example \cite{Kantorovich,OrtegaRheinboldt}.

In this paper, we consider the sequence of Newton iterates from a
quadrature perspective. The successive iterates determine a sequence of
intervals between the initial point and the root, and hence give rise to a
natural partition. This provides a direct connection
between the geometry of the Newton iteration and the mesh of the induced
quadrature rule. We first establish the basic properties of this ``Newton-generated
partition'' and then use this partition to
define the associated composite trapezoidal rule and derive its local
error. These estimates are subsequently used to study the
global behavior of the quadrature sequence and its limiting error. Later, we also obtain a more explicit bound in terms of the initial Newton decrement under a suitable smallness condition, and show that the resulting quadrature error is smaller than that of the trapezoidal rule on the original interval.
\section{Newton-generated partitions}
We begin by defining a class of ``well-behaved'' convex functions. Let $\Omega=[a,b]$ with $a<b$ and $\mathscr C^2_{m,M}(\Omega)$ denote the collection of strictly increasing functions
$f:\Omega\rightarrow\mathbb R$ satisfying:
\[
f(a)=0,\qquad
f\in C^2(\Omega),\qquad
f'(a)>0,\qquad
0<m\le f''(x)\le M
\]
for every $x\in \Omega$. For $f\in\mathscr C^2_{m,M}(\Omega)$, initialize Newton's method at $x_0=b$ and
define
\[
x_{n+1}=x_n-\frac{f(x_n)}{f'(x_n)}.
\]
Write
\[
e_n=x_n-a,\qquad
h_n=x_n-x_{n+1}
      =\frac{f(x_n)}{f'(x_n)}.
\]
The quantity $h_n$ will be referred to as the Newton decrement.
\begin{lemma}\label{lem:partition}
The sequence $(x_n)$ is bounded and strictly monotone, and
\[
[a,b]=\{a\}\cup\bigcup_{n=0}^{\infty}(x_{n+1},x_n].
\]
Consequently, the Newton iterates naturally generate a partition of the interval whose mesh widths are precisely the Newton decrements.
\end{lemma}

\begin{proof}
Since $f(x_n)>0$ and $f$ is strictly increasing, each Newton step satisfies
\[
x_{n+1}=x_n-\frac{f(x_n)}{f'(x_n)}<x_n.
\]
Taylor's theorem about $x_n$ gives
\[
0=f(a)=f(x_n)+f'(x_n)(a-x_n)
+\frac12f''(\xi_n)(a-x_n)^2,
\]
for some $\xi_n\in(a,x_n)$. Writing $e_n=x_n-a$,
\[
0=f(x_n)-f'(x_n)e_n+\frac12f''(\xi_n)e_n^2,
\]
so
\[
f'(x_n)e_n
=f(x_n)+\frac12f''(\xi_n)e_n^2>f(x_n).
\]
Hence
\[
\frac{f(x_n)}{f'(x_n)}<e_n.
\]
Therefore
\[
x_{n+1}
=x_n-\frac{f(x_n)}{f'(x_n)}
=a+e_n-\frac{f(x_n)}{f'(x_n)}>a,
\]
so $e_{n+1}>0$ and consequently $a<x_{n+1}<x_n$. Since $(x_n)$ is decreasing and bounded below by $a$, using the Monotone Convergence Theorem we get that $x_n\to L$ for some $L\ge a$. If $L>a$, then $f'(L)>0$, and passing to the limit in
\[
x_{n+1}=x_n-\frac{f(x_n)}{f'(x_n)}
\]
gives $f(L)=0$, which is impossible since $f$ is strictly increasing and $f(a)=0$. Hence $L=a$. Thus, $x_n\to a$. Consequently, every point of $(a,b]$ belongs to exactly one of the intervals $(x_{n+1},x_n]$, completing the proof.
\end{proof}
These partitions will henceforth be referred to as Newton-generated partitions.
\pagebreak
\begin{proposition}\label{prop: 2,1}
For every $n\ge0$,
\[
e_{n+1}
=\frac{f''(\xi_n)}{2f'(x_n)}\,e_n^2,
\]
where $\xi_n\in(a,x_n)$. Since $f'$ is continuous and positive on the compact interval $\Omega$, there exists $c:=\min_{x\in\Omega}f'(x)>0$. Consequently,
\[
e_{n+1}\le\frac{M}{2c}e_n^2.
\]
Moreover,
\[
h_{n+1}\le\frac{M}{2c}h_n^2.
\]
\end{proposition}

\begin{proof}
From Taylor's theorem,
\[
0=f(a)=f(x_n)-f'(x_n)e_n+\frac12f''(\xi_n)e_n^2,
\]
so
\[
f(x_n)=f'(x_n)e_n-\frac12f''(\xi_n)e_n^2.
\]
Using the Newton iteration,
\[
e_{n+1}
=e_n-\frac{f(x_n)}{f'(x_n)}
=\frac{f''(\xi_n)}{2f'(x_n)}e_n^2,
\]
which is the stated identity. Since $f''(\xi_n)\le M$ and
$f'(x_n)>0$,
\[
e_{n+1}\le\frac{M}{2c}e_n^2.
\]
For the mesh widths, expand $f(x_{n+1})$ about $x_n$
\[
f(x_{n+1})
=f(x_n)-f'(x_n)h_n+\frac12f''(\eta_n)h_n^2.
\]
Because $h_n=f(x_n)/f'(x_n)$, the linear term cancels exactly, leaving
\[
f(x_{n+1})=\frac12f''(\eta_n)h_n^2.
\]
Dividing by $f'(x_{n+1})>0$ gives
\[
h_{n+1}
=\frac{f(x_{n+1})}{f'(x_{n+1})}
\le\frac{M}{2c}h_n^2,
\]
completing the proof.
\end{proof}


\section{Newton quadrature and local error}
In this section, we consider the trapezoidal approximation of a smooth function ($f\in C^2(\Omega)$); see, for example, \cite{Atkinson}. For each local error term, we derive a Peano kernel representation and subsequently use it to obtain estimates for the corresponding error terms as a corollary. For each Newton interval, $I_n=[x_{n+1},x_n]$, we define the local trapezoidal approximation by:
\[
T_n
=
\frac{h_n}{2}
\bigl(f(x_n)+f(x_{n+1})\bigr),
\]
and set
\[
Q_N=\sum_{n=0}^{N-1}T_n.
\]
\begin{figure}[H]
    \centering
    \includegraphics[width=0.95\textwidth]{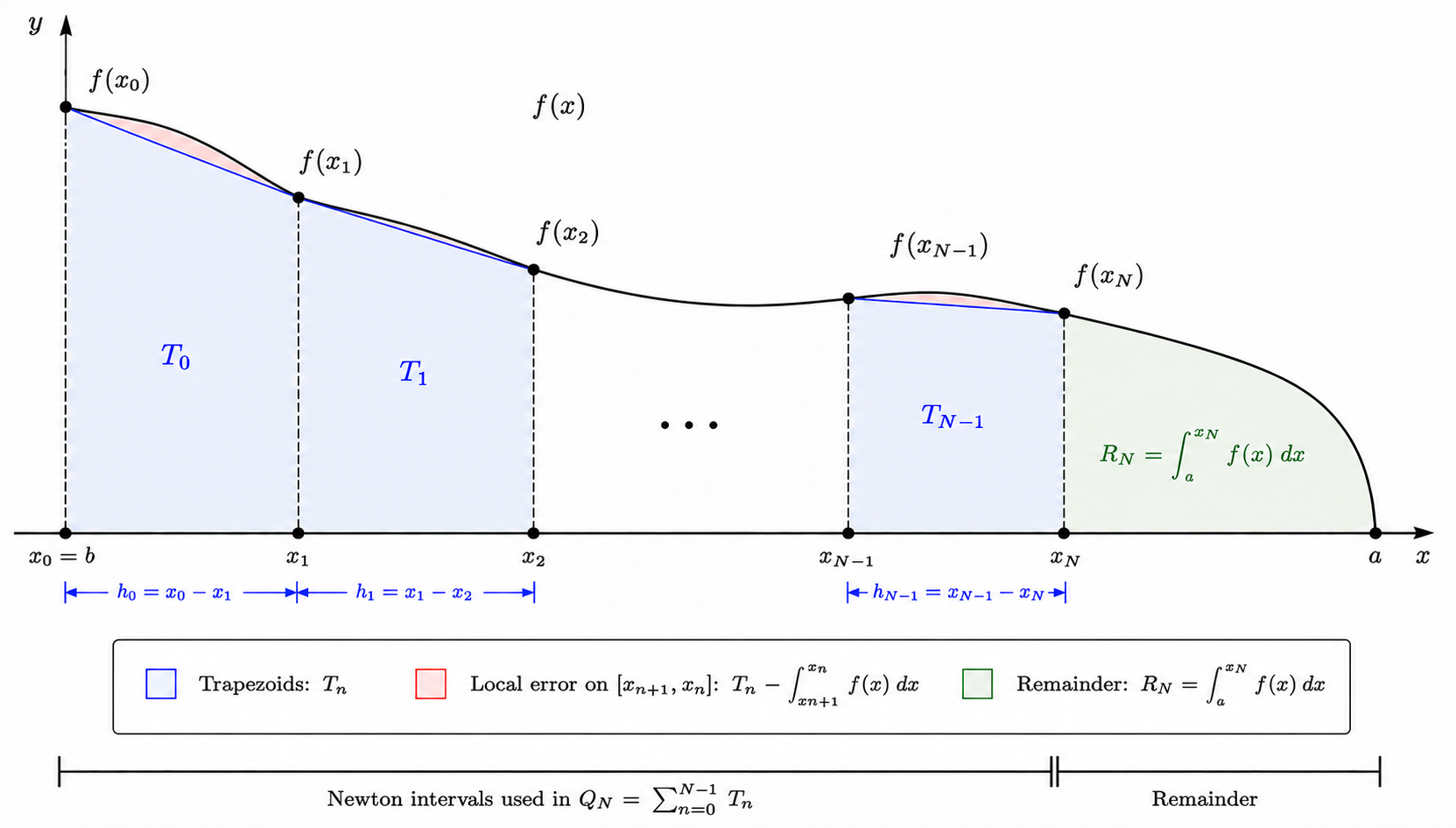}
    \caption{The Newton-generated partition and the associated composite trapezoidal approximation.}
    \label{fig:trapezoidal-decomposition}
\end{figure}
\begin{proposition}[Exact quadrature decomposition]\label{prop: 3.1}
For every $N\ge1$,
\[
\int_a^bf(x)\,dx
=
Q_N
-
\sum_{n=0}^{N-1}E_n
+
R_N,
\]
where
\[
E_n
=
T_n-\int_{x_{n+1}}^{x_n}f(x)\,dx,
\qquad
R_N=\int_a^{x_N}f(x)\,dx.
\]
\end{proposition}

\begin{proof}
Using Lemma~\ref{lem:partition},
\[
\int_a^bf(x)\,dx
=\sum_{n=0}^{N-1}\int_{x_{n+1}}^{x_n}f(x)\,dx+R_N.
\]
For each interval,
\[
\int_{x_{n+1}}^{x_n}f(x)\,dx=T_n-E_n.
\]
Summing over $n$ gives
\[
\int_a^bf(x)\,dx
=\sum_{n=0}^{N-1}T_n
-\sum_{n=0}^{N-1}E_n
+R_N
=Q_N-\sum_{n=0}^{N-1}E_n+R_N,
\]
which is exactly the desired decomposition.
\end{proof}

The local quadrature error is completely described by the Peano kernel.

\begin{proposition}[Peano-kernel representation]
For every Newton interval,
\[
E_n
=
\frac12
\int_{x_{n+1}}^{x_n}
(x-x_{n+1})(x_n-x)f''(x)\,dx.
\]
\end{proposition}

\begin{proof}
The trapezoidal error functional on an interval $[a,b]$ is
\[
L(g)=\frac{b-a}{2}\bigl(g(a)+g(b)\bigr)-\int_a^bg(x)\,dx.
\]
Since $L(1)=0$ and $L(x)=0$, the functional annihilates every affine
function. The Peano Kernel Theorem (see \cite{DavisRabinowitz}) therefore gives
\[
L(g)=\int_a^bg''(x)K(x)\,dx,
\]
where
\[
K(x)=\frac12(x-a)(b-x).
\]
Substituting $a=x_{n+1}$ and $b=x_n$ yields
\[
E_n
=\frac12\int_{x_{n+1}}^{x_n}
(x-x_{n+1})(x_n-x)f''(x)\,dx,
\]
as claimed.
\end{proof}
\begin{corollary}\label{cor: 3.3}
For every Newton interval,
\[
\frac{m}{12}h_n^3
\le
E_n
\le
\frac{M}{12}h_n^3.
\]
In particular, each local trapezoidal approximation overestimates a strictly
convex function.
\end{corollary}

\begin{proof}
From the Peano-kernel representation,
\[
E_n
=\frac12
\int_{x_{n+1}}^{x_n}
(x-x_{n+1})(x_n-x)f''(x)\,dx.
\]
Since $m\le f''(x)\le M$,
\[
\frac m2\int_{x_{n+1}}^{x_n}
(x-x_{n+1})(x_n-x)\,dx
\le E_n\le
\frac M2\int_{x_{n+1}}^{x_n}
(x-x_{n+1})(x_n-x)\,dx.
\]
The elementary integral
\[
\int_a^b(x-a)(b-x)\,dx=\frac{(b-a)^3}{6}
\]
equals $h_n^3/6$ on the interval $[x_{n+1},x_n]$, so
\[
\frac{m}{12}h_n^3
\le E_n
\le\frac{M}{12}h_n^3.
\]
Since the kernel is nonnegative and $f''>0$, every local trapezoidal
approximation overestimates a strictly convex function.
\end{proof}


\section{Global convergence and limiting error}

The cubic dependence of the local error on the Newton decrement leads to a
general convergence criterion.

\begin{theorem}\label{thm: 4.1}
Assume
\[
\sum_{n=0}^{\infty} h_n^3 < \infty.
\]
Then
\[
\lim_{N\to\infty} Q_N
=
\int_a^b f(x)\,dx
+
\sum_{n=0}^{\infty}E_n.
\]
In particular,
\[
\lim_{N\to\infty}
\left(
Q_N-\int_a^b f(x)\,dx
\right)
=
\sum_{n=0}^{\infty}E_n
>0.
\]
\end{theorem}

\begin{proof}
The local estimate gives
\[
|E_n|
\leq \frac{M}{12}h_n^3.
\]
Since
\[
\sum_{n=0}^{\infty} h_n^3 < \infty,
\]
the comparison test implies
\[
\sum_{n=0}^{\infty}|E_n|
\leq
\frac{M}{12}\sum_{n=0}^{\infty}h_n^3
<\infty.
\]
Thus the error series converges absolutely, and hence
\[
\lim_{N\to\infty}
\sum_{n=0}^{N-1}E_n
=
\sum_{n=0}^{\infty}E_n.
\]
Moreover,
\[
\lim_{N\to\infty}x_N=a.
\]
Since \(f\) is continuous,
\[
\lim_{N\to\infty}R_N
=
\lim_{N\to\infty}\int_a^{x_N}f(x)\,dx
=
\int_a^a f(x)\,dx
=
0.
\]
Taking the limit as \(N\to\infty\) in the exact decomposition
\[
\int_a^b f(x)\,dx
=
Q_N-\sum_{n=0}^{N-1}E_n+R_N
\]
yields
\[
\lim_{N\to\infty}Q_N
=
\int_a^b f(x)\,dx
+
\sum_{n=0}^{\infty}E_n,
\]
as required.
\end{proof}
\begin{corollary}
Under the standing assumptions on $f$, the sequence $(Q_N)$ converges and
\[
\lim_{N\to\infty} Q_N
=
\int_a^b f(x)\,dx + E_\infty,
\qquad
E_\infty:=\sum_{n=0}^{\infty}E_n.
\]
Moreover, there exist constants $C>0$ and $N_0\in\mathbb{N}$ such that, for
all $N\geq N_0$,
\[
\left|
Q_N-
\left(
\int_a^b f(x)\,dx+E_\infty
\right)
\right|
\leq C e_N^2.
\]
\end{corollary}

\begin{proof}
By Theorem~\ref{thm: 4.1},
the series $\sum_{n=0}^{\infty}E_n$ converges and
\[
\lim_{N\to\infty}Q_N
=
\int_a^b f(x)\,dx+E_\infty.
\]
From the exact decomposition in Proposition~\ref{prop: 3.1},
\[
Q_N-
\left(
\int_a^b f(x)\,dx+E_\infty
\right)
=
R_N-\sum_{n=N}^{\infty}E_n.
\]

Let
\[
L:=\max_{x\in[a,b]}|f'(x)|.
\]
Since $f(a)=0$, the mean value theorem yields
\[
|f(x)|
=
|f(x)-f(a)|
\leq L(x-a),
\qquad x\in[a,b].
\]
Hence
\[
|R_N|
=
\left|\int_a^{x_N}f(x)\,dx\right|
\leq
\int_a^{x_N}L(x-a)\,dx
=
\frac{L}{2}e_N^2.
\]

On the other hand, Proposition~\ref{prop: 2,1} gives
\[
h_{n+1}\leq \frac{M}{2c}h_n^2.
\]
Since $h_n\to0$, there exists $N_0\in\mathbb{N}$ such that
\[
\frac{M}{2c}h_n\leq\frac12,
\qquad n\geq N_0.
\]
Thus
\[
h_{n+1}\leq\frac12h_n,
\qquad n\geq N_0.
\]
It follows that, for $N\geq N_0$,
\[
\sum_{n=N}^{\infty}h_n^3
\leq
h_N^3\sum_{k=0}^{\infty}\frac{1}{2^{3k}}
=
\frac{8}{7}h_N^3.
\]
By Corollary~\ref{cor: 3.3},
\[
\sum_{n=N}^{\infty}E_n
\leq
\frac{M}{12}\sum_{n=N}^{\infty}h_n^3
\leq
\frac{2M}{21}h_N^3.
\]
Since
\[
h_N=x_N-x_{N+1}\leq x_N-a=e_N,
\]
we obtain
\[
\sum_{n=N}^{\infty}E_n
\leq
\frac{2M}{21}e_N^3.
\]
Therefore, for $N\geq N_0$,
\[
\begin{aligned}
\left|
Q_N-
\left(
\int_a^b f(x)\,dx+E_\infty
\right)
\right|
&\leq
|R_N|+\sum_{n=N}^{\infty}E_n\\
&\leq
\frac{L}{2}e_N^2
+
\frac{2M}{21}e_N^3.
\end{aligned}
\]
Since $e_N\leq b-a$,
\[
\left|
Q_N-
\left(
\int_a^b f(x)\,dx+E_\infty
\right)
\right|
\leq
\left(
\frac{L}{2}
+
\frac{2M}{21}(b-a)
\right)e_N^2.
\]
Thus the asserted estimate holds with
\[
C=
\frac{L}{2}
+
\frac{2M}{21}(b-a).
\]
This completes the proof.
\end{proof}
The preceding results describe the convergence of the Newton-generated
trapezoidal rule. We now examine a different phenomenon: the contribution of
the initial Newton interval, which is never subsequently subdivided.

\begin{proposition}
The local error on the first Newton interval satisfies
\[
E_0
\ge
\frac{m}{12}
\left(\frac{f(b)}{f'(b)}\right)^3.
\]
Consequently, the initial Newton interval contributes a strictly positive cubic error bounded below in terms of the first Newton decrement.
\end{proposition}

\begin{proof}
From Corollary~\ref{cor: 3.3},
\[
E_0\ge\frac{m}{12}h_0^3.
\]
Since
\[
h_0=x_0-x_1=\frac{f(b)}{f'(b)},
\]
substituting this identity gives
\[
E_0
\ge
\frac{m}{12}
\left(\frac{f(b)}{f'(b)}\right)^3,
\]
which proves the claim.
\end{proof}

This contribution persists because the interval $[x_1,b]$ is created during the
first Newton step and is never subsequently subdivided by the Newton iteration.



\section{Error Analysis Under Small Newton Decrement}
The results of the previous sections show that the Newton-generated trapezoidal rule has a finite limiting error, although this error need not vanish as the number of Newton steps tends to infinity. We now give an explicit bound for this accumulated error and identify a condition under which it is smaller than the error of the one-panel trapezoidal rule on the original interval. The quantity
\[
\rho:=\frac{h_0}{b-a}
=\frac{f(b)}{(b-a)f'(b)}
\]
will be used to measure the size of the initial Newton decrement relative to the length of the interval. In this section, we assume the following quantity is positive, i.e.
\[
c:=\min_{x\in[a,b]}f'(x)>0.
\]

\begin{theorem}\label{thm:5.1}
Suppose
\[
\frac{Mh_0}{2c}<1.
\]
Then the limiting error
\[
E_\infty:=\sum_{n=0}^{\infty}E_n
\]
satisfies
\[
0<E_\infty
\leq
\frac{Mh_0^3}
{12\left(1-\left(\frac{Mh_0}{2c}\right)^3\right)}.
\]
Equivalently,
\[
0<E_\infty
\leq
\frac{M(b-a)^3\rho^3}
{12\left(1-\left(\frac{M(b-a)\rho}{2c}\right)^3\right)}.
\]
\end{theorem}

\begin{proof}
By Proposition~\ref{prop: 2,1},
\[
h_{n+1}\leq\frac{M}{2c}h_n^2.
\]
Set
\[
q:=\frac{Mh_0}{2c}<1.
\]
An induction gives
\[
h_n\leq h_0q^{2^n-1}\leq h_0q^n.
\]
Hence
\[
\sum_{n=0}^{\infty}h_n^3
\leq
h_0^3\sum_{n=0}^{\infty}q^{3n}
=
\frac{h_0^3}{1-q^3}.
\]
Since, by Corollary~\ref{cor: 3.3},
\[
0<E_n\leq\frac{M}{12}h_n^3,
\]
we obtain
\[
0<E_\infty
\leq
\frac{M}{12}\sum_{n=0}^{\infty}h_n^3
\leq
\frac{Mh_0^3}{12(1-q^3)}.
\]
Substituting the definition of $q$ proves the result.
\end{proof}
\flushleft We next compare this bound with the error of the one-panel trapezoidal rule
\[
T[f]:=\frac{b-a}{2}\bigl(f(a)+f(b)\bigr).
\]
Let
\[
E_{\mathrm{unif}}
:=
T[f]-\int_a^b f(x)\,dx.
\]
Since $f''(x)\geq m$ on $[a,b]$,
\[
E_{\mathrm{unif}}\geq\frac{m(b-a)^3}{12}.
\]

\begin{corollary}
Suppose
\[
\frac{Mh_0}{2c}<1
\]
and
\[
\frac{M}{m}
\frac{\rho^3}
{1-\left(\frac{M(b-a)\rho}{2c}\right)^3}
<1.
\]
Then
\[
E_\infty<E_{\mathrm{unif}}.
\]
\end{corollary}

\begin{proof}
By Theorem~\ref{thm:5.1},
\[
E_\infty
\leq
\frac{M(b-a)^3}{12}
\frac{\rho^3}
{1-\left(\frac{M(b-a)\rho}{2c}\right)^3}
<
\frac{m(b-a)^3}{12}
\leq E_{\mathrm{unif}}.
\]\end{proof}

\begin{remark}
The quantity $\rho$ measures the relative size of the first Newton step. Thus, the above condition shows that when the initial Newton decrement is sufficiently small and the curvature is suitably controlled, the total bias of the Newton-generated quadrature can be smaller than the error of the one-panel trapezoidal rule on the original interval.
\end{remark}

\section{Conclusion}

We have shown that the quadrature induced by the Newton iterates converges
to a limit with a strictly positive bias. The local error is controlled by
the Newton decrements, and under a suitable smallness condition on the
initial decrement, the limiting bias admits an explicit bound. In
particular, this bound yields a condition under which the limiting bias is
smaller than the error of the one-panel trapezoidal rule on the original
interval. These results illustrate how the convergence properties of
Newton's method are reflected in the quadrature rule generated by its
iterates.
\vspace{0.25cm}
\paragraph{\textbf{AI Disclosure:}}
This paper was originally written in 2022, when the author
was still a high school student. Except for the generation of Figure~1, no
artificial intelligence tools were used in the writing, revision, or preparation of this paper.

\end{document}